\documentclass[]{article}       

\usepackage{graphicx}
\usepackage{natbib}
\usepackage{amsmath} 
\usepackage{amsfonts}
\usepackage{hyperref}
\usepackage{caption}
\usepackage{subcaption}

\usepackage{standalone}
\usepackage{tikz}
\usetikzlibrary{shapes,arrows}
\usepackage{booktabs}

\begin{document}

\title{Learning to Handle Parameter Perturbations in Combinatorial Optimization: an Application to Facility Location}

\author{Andrea Lodi$^1$ \and Luca Mossina$^2$ \and Emmanuel Rachelson$^2$}
\date{$^1$ CERC and \'Ecole Polytechnique de Montr\'eal,  \texttt{andrea.lodi@polymtl.ca}  \\
      $^2$ ISAE-SUPAERO, Universit\'{e} de Toulouse,  \texttt{name.surname@isae-supaero.fr}
}

\maketitle

\begin{abstract}
We present an approach to couple the resolution of Combinatorial Optimization problems with methods from Machine Learning, applied to the single source, capacitated, facility location problem.
Our study is framed in the context where a reference facility location optimization problem is given.
Assuming there exist data for many variations of the reference problem (historical or simulated) along with their optimal solution, we study how one can exploit these to make predictions about an unseen new instance.
We demonstrate how a classifier can be built from these data to determine whether the solution to the reference problem still applies to a new instance.
In case the reference solution is partially applicable, we build a regressor indicating the magnitude of the expected change, and conversely how much of it can be kept for the new instance.
This insight, derived from a priori information, is expressed via an additional constraint in the original mathematical programming formulation.
We present an empirical evaluation and discuss the benefits, drawbacks and perspectives of such an approach.
Although presented through the application to the facility location problem, the approach developed here is general and explores a new perspective on the exploitation of past experience in combinatorial optimization.
\end{abstract}

\section{Introduction}
Solving combinatorial optimization problems can be highly challenging.
Many problems have been shown to be computationally hard to solve \citep{Garey1979ComputersAI}, that is, no polynomial-time algorithms exist.
Despite this intrinsic complexity, the state-of-the-art modeling and solving solutions can handle real-world instances via exact methods, heuristics or a mix of the two.

Our rationale goes as follows: given a combinatorial problem, we assume the existence of a reference problem \textit{instance} $P_{ref}$, derived from the description of a system in its nominal state.
This reference can be modeled via Mixed Integer Programming (MIP) and solved to its optimum $x_{ref}^*$ via state-of-the-art solvers based on branch-and-bound (BB) methods like CPLEX \citep{cplex_solver}.
If a modification of the reference parameters occurs, the problem needs to be solved again.
We operate as if we had to respond to this parameter perturbation (also referred to as \textit{disruption}) under a tight optimization time budget, aiming at accelerating the descent towards the optimum.
Given the modifications that have occurred in the past or that can be simulated \textit{a priori}, one can build a dataset of resolutions.
Through the example of the Facility Location Problem, we study how one can extract information from this dataset to facilitate the resolution of future instances.
We cast this problem as a Supervised Learning problem to predict whether a modification will affect and to what extent the reference optimal solution.
According to this prediction, additional constraints are added to the MIP formulation of the new instance. 
In order to evaluate the approach, the two formulations are fed to the solver and the performance are compared. 

Section \ref{sec:l4opt-rev} reviews the related literature on applying Machine Learning to Optimization problems.
Section \ref{sec:FLP} provides a description of the Facility Location Problem and Section \ref{sec:learnbound} introduces how the prediction problem can be formulated mathematically and linked to existing Machine Learning approaches.
Experiments and numerical results are presented and discussed in Section \ref{sec:methods}.
Section \ref{sec:conclusion} discusses the advantages and drawbacks of the proposed method and present future research directions.

\section{Machine Learning for Optimization}\label{sec:l4opt-rev}
The remarkable results achieved in recent years by Machine Learning (ML), for example in the field of image recognition \citep{krizhevsky2012imagenet} or Reinforcement Learning \citep{mnih2015human}, sparked a lot of interest in the Operations Research community.
Some results have emerged in applications such as guiding the branching process in branch-and-bound optimization algorithms, in the form of node exploration ordering \citep{he2014learning} and approximating performant branching rules \citep{alvarez2017machine, KhaBodSonNemDil16, lodi2017learning}.
In \citet{Liberto2016DASH} one can find algorithmic ideas on how to handle the many existing heuristics already developed for branch-and-bound methods.

The traveling salesman problem has been addressed via an ad-hoc neural network architecture in a Reinforcement Learning context \citep{bello2017neural, khalil2017learning, deudon2018learning}, aiming at learning to build incrementally a solution.

Another field of application is the interaction between existing optimization solvers and \textit{a priori} knowledge on the problem.
For instance, \citet{kruber2017learning} find an approach to automatically handle MIP problem decompositions, by detecting if and what reformulation to apply within a dedicated software.
On the same topic, \citet{Basso2018} bring evidence to why such an approach can be carried out, showing empirically that the role of an expert needed when using decomposition methods can be, at least in part, automated via learning algorithms.
In \citet{bonami2018learning}, the authors present a method to select the best resolution method for a Mixed Integer Quadratic Programming problem from the different algorithms offered within the CPLEX framework.

Machine Learning has proven useful in producing a description of the optimal solution of a yet unseen instance.
\citet{fischetti2018machine} describe a real-world application with the problem of a wind turbine park layout.
Minimizing the complex turbulence induced by the positioning of wind turbines generates a series of difficult combinatorial problems.
Where different configurations would need to be evaluated at a considerable computational cost, the authors approximate the solution values to these optimizations via ML, thanks to a dataset built \textit{a priori}.
On a similar line lies the work of \citet{larsen2018predicting}.
When the problem of choosing the optimal load planning for containers on freight trains cannot be solved online because of time constraints and insufficient information (tactical level), a set of offline cases can be collected and used to get a description of the solution of a new instance, at an aggregate level.
Such description provides meaningful insights to decision makers in a real-time context. 

Strictly connected with our work, \citet{Ahmed2019} consider a number of ML techniques to extract information from previously solved instances of Unit Commitment problems and leverage such pieces of information to improve the MIP performance when solving similar instances again and again.

The interested reader is referred to \citet{BLP18} for a recent survey on the subject.

\section{Constrained Facility Location Problems} \label{sec:FLP}
Our application of choice will be that of the Capacitated Facility Location Problem (CFLP).
Given a group of customers and a list of potential facilities that can serve them, one must open a certain amount of facilities and assign a unique facility to each of the customers. The goal of the problem is to satisfy the demand of the customers while minimizing the associated operational and start-up costs.
The facilities are constrained in capacity, that is, the quantity of total service each can provide is limited.
When solving this problem one determines which facilities must be opened and which customers they will serve. 

\subsection{Mathematical Programming Formulation} \label{sec:mip-flp}
The CFLP can be written as a Binary Integer Program, with decision variables $y_{j}, x_{ij} \in \{0,1\}$, indicating respectively whether facility $j$ is activated (also referred to as \textit{open}) and whether customer $i$ is served by $j$.
\begin{align}
            & \underset{x, y}{\text{minimize}} & & \sum_{i \in I} \sum_{j \in J} c_{ij} x_{ij} + \sum_{j \in J} f_j y_j & \label{eq:flp1}\\
            & \text{subject to}             & & \sum_{i \in I} d_{i} x_{ij} \leq s_j y_j  & \forall j \in J        \label{eq:flp2}\\
            &                               & & \sum_{j \in J} x_{ij} = 1                 & \forall i \in I \label{eq:flp3}\\
            &                               & & x_{ij} \in \{0,1\}, \ y_{j}  \in \{0,1\} & \forall i \in I, \forall j \in J \label{eq:flp4}\\
            &                               & & I = \{1, 2, \dots, N_C\}, J = \{1, 2, \dots, N_F\}, &\label{eq:flp5}
\end{align}

\noindent where $i \in I$ are the customers and $j \in J$ the available facilities;
$c_{ij} \geq 0$ is the cost of serving customer $i$ from $j$;
$f_j \geq 0$ is the fixed cost for opening facility $j$;
$d_i \geq 0$ is the demand of customer $i$;
$s_j \geq 0$ is the capacity of facility $j$;
$N_F$ and $N_C$ are the number of facilities and
of customers, respectively.

The objective function (\ref{eq:flp1}) quantifies the total cost of a given assignment, composed of the fixed start-up costs $\sum_{j \in J} f_j y_j $ and $\sum_{i \in I} \sum_{j \in J} c_{ij} x_{ij}$, the operational costs.
Constraints (\ref{eq:flp2}) ensure that the total demand of the customers assigned to $j$ does not exceed its capacity $s_j$.
Constraints (\ref{eq:flp3}) ensure that each customer is served by exactly one facility.
For convenience, we write an instance as $P_{}=\left(\boldsymbol{c}_{},\boldsymbol{f}_{},\boldsymbol{d}_{},\boldsymbol{s}_{}\right)$.

In the various forms of the Facility Location Problem (FLP), the binary variables associated with facility activation and customer-to-facility assignment make the problem computationally challenging.
Even the uncapacitated version, for instance, has been proven to be NP-hard \citep{cornuejols1983uncapacitated}. 
For a comprehensive treatment on exact formulations of the FLP, we refer the reader to the work of \citet{klose2005facility}.
Beyond mathematical programming and exact methods, heuristics for the FLP have received a fair amount of attention \citep{cornuejols1983uncapacitated,guastaroba2012kernel,GUASTAROBA2014_heuristic}.

\subsection{Perturbations in a Reference Problem}
We assume the existence of a reference instance $P_{ref}=\left(\boldsymbol{c}_{ref},\boldsymbol{f}_{ref},\boldsymbol{d}_{ref},\boldsymbol{s}_{ref}\right)$ where the parameters are considered to be in their nominal state. 
This problem has $\{\boldsymbol{x}_{ref}^{*}, \boldsymbol{y}_{ref}^*\}$ as its reference solution.
In case a disruption in the system occurs, for example an anomaly in the capacity, we wish to learn whether the reference solution, or part of it, is still applicable without running the full optimization.
To that end, we wish to predict what parts of the solution need to be changed by providing hints to the solver, via an additional constraint.
We shall write a new instance $P'=\left(\boldsymbol{c}',\boldsymbol{f}',\boldsymbol{d}',\boldsymbol{s}'\right)$, which we see as a variation $P' = P_{ref}+\Delta P$, where a disruption $\Delta P$ affected one or more of $P_{ref}$'s parameters.
In practice, one would proceed as in the following example.
The given (simple) $P_{ref}$  with capacities $\boldsymbol{s}_{ref} = (5000, 5000, 5000)$ is affected by a disruption and the derived $P'$ has capacities $\boldsymbol{s}' = (5000, 5088, 4304)$.
The ML algorithm will take $\Delta P = \{0, 88, -696\}$ as input and predict how $\{\boldsymbol{x'}_{}^{*}, \boldsymbol{y'}_{}^*\}$ could differ from $\{\boldsymbol{x}_{ref}^{*}, \boldsymbol{y}_{ref}^*\}$, before running the optimization for $P'$ (see Section \ref{sec:learnbound}).
In our case, the hint to the solver will specify how many (if any at all) of the facilities active in $P_{ref}$ will be still active, without specifying which of the three.

The FLP can be adapted to different domains.
Consider for instance the airport system of a region.
If a runway at a major airport has problems or needs maintenance, one could think of this as a disruption.
One could then want to learn and predict how to react to the closure of a runway, how much of the traffic would need to be delayed, rescheduled or rerouted or, in general, how to be guided to manage the disruption.

\section{Learning constraints}\label{sec:learnbound}
In this Section, we formalize the search for a constraint that will be used to accelerate the resolution of FLP instances facing disruptions. We formulate this as a statistical learning problem and describe how we couple it with MIP resolution.

\subsection{Learning to Bound the Changes to a Reference Solution}
Given the information gained from past resolutions or simulated data, we  investigate the possibility of predicting the number of facilities we expect to change.
We want to predict that ``only a proportion $\gamma=0.15$ of the facilities needs replanning", and impose that as a constraint. Indeed, it is well known in the case of MIPs with binary variables \citep{localbranching} that it is easy  to express the neighborhood of a feasible solution by a linear constraint and the resulting MIP is generally way easier to solve than the original problem. In addition, for FLP it is also well known that the critical choice is associated with the facilities to open while the assignment of clients to facilities is easier to deal with.

Then, the problem boils down to predicting the proportion of facilities opened in $\boldsymbol{y}_{ref}^*$ to be kept open in $\boldsymbol{y'}^*$, the optimal solution of the perturbed instance $P'$.
This estimation is a scalar value $\widehat{\gamma}$ determined via a function $\Gamma(\Delta P)$,  where $\Gamma : \mathbb{R} \rightarrow [0,1]$.
Given the $\left\{\left(\Delta P_k, \gamma_k\right)\right\}_{1\leq k\leq K}$ data from the $K$ past resolutions, this is a regression problem, which we can tackle with any Statistical Learning method \citep{hastie2009unsupervised}.

The drawback of such an approach is the potential introduction of bias in the solution of $P'$ if the true optimum of $P'$ is far from  $\{\boldsymbol{x}_{ref}^{*}, \boldsymbol{y}_{ref}^*\}$ and the ML fails to detect it or if 
it cuts the optimum from the admissible domain.
In fact, in case the parameters $\left(\boldsymbol{c}',\boldsymbol{f}',\boldsymbol{d}',\boldsymbol{s}'\right)$ generating $P'$ are very different from $P_{ref}$, a good 
model would impose no additional constraint, thus falling back to the full problem.

The information on $\widehat{\gamma} = \Gamma(\Delta P)$ can be modeled through the linear constraint (\ref{eq:constr}), thus yielding the new problem
\begin{align}
\nonumber  & \underset{x,y}{\text{minimize}} & & \sum_{i \in I} \sum_{j \in J} c_{ij} x_{ij} + \sum_{j \in J} f_j y_j & \\
\nonumber  & \text{subject to}             & & \textrm{all previous constraints and} \\ 
        &                               & & \sum_{r \in R_{ref}} y_{r} \geq \widehat{\gamma} \    N_{ref}, &
  \label{eq:constr}
\end{align}

\noindent where $R_{ref} = \{j \in J \mid y_{ref,j}^* = 1\}$ is the set of indices of the facilities opened in the reference solution, $\widehat{\gamma} \in [0,1]$ is the parameter predicted via ML and $N_{ref} = |R_{ref}|$ is the number of facilities activated in the reference problem.

We do not attempt to identify directly which of the facilities need to change or stay as in $\{\boldsymbol{x}_{ref}^{*}, \boldsymbol{y}_{ref}^*\}$, but rather to restrict the number of facilities that can be changed.
If $\widehat{\gamma} \approx 0$, the disruption data conveys no information as to the new optimum and all the facility allocations should be left to the solver.
In that case, it is likely that the variation on $P_{ref}$ is so strong that prior information is of no help and we need to solve the new instance as a full problem in its original formulation.
Conversely, if $\widehat{\gamma} \approx 1$, we conclude that all of the reference facilities are to be left activated, without excluding further facilities from being activated.

\subsection{Structuring the prediction problem} \label{how_pred}
The overall problem of predicting  $\widehat{\gamma} = \Gamma(\Delta P)$ is a regression problem.
Previously solved instances $P_k$ provide values $(\Delta P_k, \gamma_k)$ to assemble a training set.
Then, most ML methods will search for $\Gamma$ by minimizing a loss function defined as an expected value of individual losses over the training set, such as the least squares.
Such methods can thus be sensitive to the training samples' distribution and might overfit to the majority value.
In the specific case of predicting the proportion of facilities to keep open, the distribution of observed values for samples $(\Delta P_k,\gamma_k)$ is strongly affected by the fact that in many cases, \emph{all} the facilities in the reference solution should remain open (for details, see Section \ref{sec:gendata}).
To compensate for this imbalance, our predictor's architecture features two levels.
\begin{figure}
	\centering
    \begin{tabular}{ll}
        \begin{minipage}{4cm}
        \includegraphics[width=0.8\linewidth]{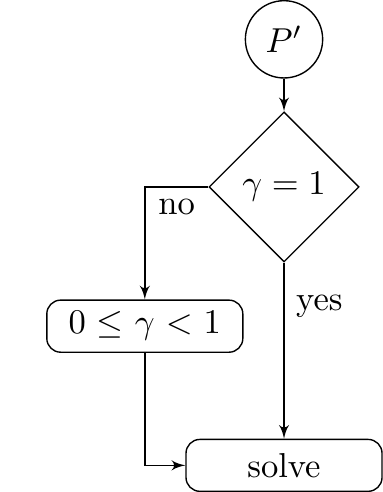}
        \end{minipage}   &
        \begin{minipage}{7cm}
		\textbf{Learned bound}
		\begin{enumerate}
		\item \texttt{New instance} $P'$.
		\item \texttt{if} $\widehat{\gamma} = \Gamma(P') = 1$:
		\item []\ \ \ \texttt{add} $y_{r}  = 1,  \forall r \in R_{ref}$ 
		\item  \texttt{else if} $\widehat{\gamma} = \Gamma(P') \in \left[0,1\right)$:
		\item []\ \ \ \texttt{add}$\sum_{r \in R_{ref}} y_{r} \geq \widehat{\gamma} \    N_{ref}$
		\item \texttt{Solve}
		\end{enumerate}
        \end{minipage}
    \end{tabular}

    \caption{Learning to bound the FLP}
    \label{fig:diagram}
\end{figure}

First, we fit a binary classifier to indicate whether the open facilities in $\{\boldsymbol{x}_{ref}^{*}, \boldsymbol{y}_{ref}^*\}$ should all be kept open.
If this is the case, then $\widehat{\gamma} := 1$ without any further computation.
Second, if we cannot assert that all the reference active facilities need to be left active, we call a specific regression function.
This $\Gamma(\Delta P)$ function is trained on all the data points except for those for which $\gamma_k=1$.
Then, given a new problem $P'$, we predict how much of the reference solution has to be changed, adding that as constraint (\ref{eq:constr}) to the original problem.
Figure \ref{fig:diagram} summarizes the $\widehat{\gamma}$ estimation process and the resolution of a new instance $P'$.

\section{Experimental Evaluation} \label{sec:methods}
Given the reference instance $P_{ref}$ and an instance $P'$ derived by perturbing $P_{ref}$, our aim is to provide a good solution to $P'$ within a small time budget.

Our $P_{ref}$ is taken from the CFLP instances of \citet{GUASTAROBA2014_heuristic}\footnote{Data available at \url{http://or-brescia.unibs.it/instances/instances_sscflp}}.
We test our method on three different reference cases, respectively $capa_1$, $capb_1$ and $capc_1$.
For conciseness, we refer to these as $A, B, C$, and their corresponding MIP models as $P_{ref}^A, P_{ref}^B, P_{ref}^C$.
In each case, we have $N_F=100$ facilities available that share the same capacity $s^{A} = s_1^{A} = s_2^{A} = \dots = s_{100}^{A}$.

\subsection{Generation of Training Data}\label{sec:gendata}
For each of the three instances, we produce three datasets of sizes $N_A = N_B = N_C =10000$.
Given $P_{ref}$, we generate our training data by perturbing the capacities of the facilities.
We randomly perturb between 5\% and 50\% of the facilities as follows.

\begin{enumerate}
    \item Pick uniformly at random between 5\% and 50\% of the $N_F=100$ facilities.
    \item For each facility $j$, apply a Gaussian noise to its capacity: \\\quad $s_j \leftarrow s_j + \mathcal{N}(0, \sigma=0.2 \times s_j)$.
    \item All the remaining data is left unchanged. This yields $\Delta P_k$.
    \item Solve the new instance $P_k = P_{ref}+\Delta P_k$ with a MIP solver and get the number $Y_k$ of facilities that were in $x_{ref}^*$ and are also in $x_{k}^*$.
    \item Add the new point $(\Delta P_k, \gamma_k = Y_k / N_F)$ to the dataset.
\end{enumerate}

After generating the data, we can observe the distribution of $Y$, the number of original open facilities still open in the optimal solution of the randomized instances.
At their nominal values, the optimal solutions for $A$, $B$ and $C$ have respectively 7, 11 and 11 facilities open.
For instance, given the $k$-th point in our generated data from case $A$, if $Y_k=4$, this means that after the perturbation of the capacities, of the seven facilities open in the reference optimal solution, four were included in the optimal solution of the derived problem $P_k$.

\begin{figure}
  \centering
    \begin{subfigure}{.35\textwidth}
      \centering
      \includegraphics[width=1\linewidth]{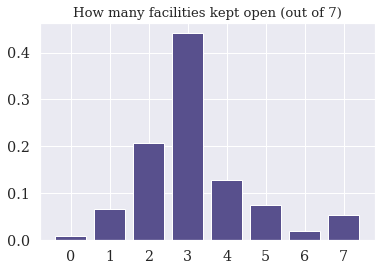}
      \caption{Reference $A$}
      \label{fig:hist_a_sub2}
    \end{subfigure}
    \begin{subfigure}{.35\textwidth}
      \centering
      \includegraphics[width=1\linewidth]{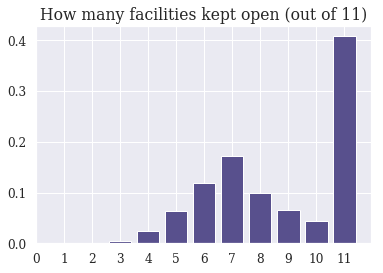}
      \caption{Reference $B$}
      \label{fig:hist_b_sub2}
    \end{subfigure}

\centering
\begin{subfigure}{.35\textwidth}
  \centering
  \includegraphics[width=1\linewidth]{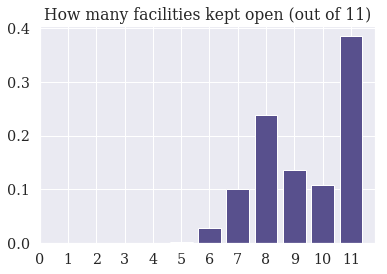}
  \caption{Reference $C$}
  \label{fig:hist_c_sub2}
\end{subfigure}

\caption{Number of facilities kept after reoptimization}
\label{fig:hist_Y}
\end{figure}

Figure \ref{fig:hist_Y} presents three different scenarios.
With instance $A$ we see that on average less than half of the 7 original facilities are kept open, that is, the perturbations affect the new optimal solutions.
In the case of $B$, for most of the data points the original reference solution is still optimal after the perturbations, as seen in the right-most bar of Figure \ref{fig:hist_b_sub2}.
Case $C$ is intermediate, with at least eight facilities out of eleven still open after the variations of the capacities.

\subsection{Learning the Prediction Model}
As mentioned in Section \ref{how_pred}, we need a binary classifier and a regressor to account for sample distribution imbalance.
The number of data for binary classification depends on the shape of the generated data, that is, on how many extreme cases we observe (rightmost column in histograms of Figure \ref{fig:hist_Y}).

We split the training data into classification and regression data as follows.
In case A, we have $717$ data points of value $7$, corresponding to the observations where all the facilities open in the reference solution are to be left open in the perturbed instances.
Thus, we set these $717$ points aside and randomly sampled an equal amount of points among the other points, yielding a learning data set of $N_{Classification}A = 1434$.
The remaining $9283$ points form the regression learning data.
\begin{table}
\centering
\begin{tabular}{rccc}
\textbf{Reference} & \textbf{Classification}  & \textbf{Regression}  & \textbf{Total} \\ \midrule
A        & 1434                & 9283            & 10000 \\
B        & 8056                & 5965            & 10000 \\
C        & 7686                & 6179            & 10000
\end{tabular}
\caption{Training Data: repartition between binary classification and regression}
\label{tab:partition}
\end{table}
This process was repeated for all the three cases, yielding the distributions reported in Table \ref{tab:partition}.

In accordance to standard practice in ML \citep{hastie2009unsupervised}, at learning time the data set was randomly split into a \textit{training} and a \textit{test} partition.
The predictive models were fit on the training data but evaluated on the test data.

\subsubsection{Comparison of ML algorithms: Classification}\label{sec:ml_bin}
For the binary classification task, we tested a set of binary classifiers among which Extremely Randomized Trees (Extra Trees) \citep{geurts2006extremely}, Neural Network classifier \citep{Goodfellow-et-al-2016}, Logistic Regression classification and Naive Bayes classification \citep{hastie2009unsupervised}.

\begin{table}
\centering
\begin{tabular}{rrcc} 
\textbf{Reference}   &\textbf{Model}       & \multicolumn{1}{c}{\textbf{Accuracy}} & \multicolumn{1}{c}{\textbf{False Positive Errors}} \\
 \midrule
A  &Extra Trees & 0.756 & 0.106                                   \\
   &Neural Net  & 0.717 & 0.160                                   \\ 
   &Logistic    & 0.648 & 0.189                                   \\ 
   &Naive Bayes & 0.641 & 0.191                                  \\   
   \midrule                                                 
B&Extra Trees & 0.840 & 0.049 \\   
&Neural Net  & 0.857 & 0.083 \\
&Naive Bayes & 0.782 & 0.100 \\
&Logistic    & 0.629 & 0.192       \\ 
    \midrule                            
C & Extra Trees & 0.864 & 0.033                                   \\
  & Neural Net  & 0.858 & 0.085                                   \\
  & Naive Bayes & 0.776 & 0.104                                   \\
  & Logistic    & 0.616 & 0.199                                  \\
\end{tabular}
\caption{Comparison of binary classifiers}
\label{tab:classif_comp}
\end{table}

As can be seen across Table \ref{tab:classif_comp}, Extra Trees emerged as the best performing in terms of accuracy and false positives errors, that is, over the total of the prediction, how many were predicted as false positives.
Here, false positives are cases where the reference solution's open facilities should \textit{not} be all kept open in the new instance but the classifier prescribes otherwise.
This could introduce bias, making this metric important for our task.

\subsubsection{Comparison of ML algorithms: Regression}\label{sec:ml_reg}
For the regression task we selected Extremely Randomized Trees out of a pool of models comprising Multiple Linear Regression, Extremely Randomized Trees and a Multi-layer Perceptron artificial Neural Network.

\begin{table}
\centering
\begin{tabular}{rccc}
\textbf{Reference}      & \textbf{Linear Regression}    & \textbf{Extra Trees}   & \textbf{Neural Net}    \\ \midrule
$A$ & 0.998 & {0.883} & 0.975  \\
$B$ & 2.087 & {1.780} & 2.099  \\
$C$ & 1.114 & {0.750} & 1.097   
\end{tabular}
\caption{Comparison of regression models by Mean Squared Error}
\label{tab:reg_mse}
\end{table}

In Table \ref{tab:reg_mse} are the values of the Mean Squared Error (MSE) computed on the test partition.
Extra Trees yielded the best performance of the tested methods.
While Neural Networks have recently been the object of intense research and remarkable results, we think its poorer results can explained by the limited amount of training data ($N \leq 10000$) and the highly nonlinear relations between features and response variable, for which more data could have been needed.

\begin{figure}
\centering
\begin{subfigure}{.4\textwidth}
  \centering
  \includegraphics[width=1\linewidth]{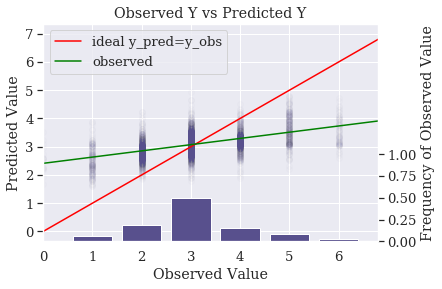}
  \caption{}
  \label{fig:scatter_a_sub1}
\end{subfigure}%
\begin{subfigure}{.4\textwidth}
  \centering
  \includegraphics[width=1\linewidth]{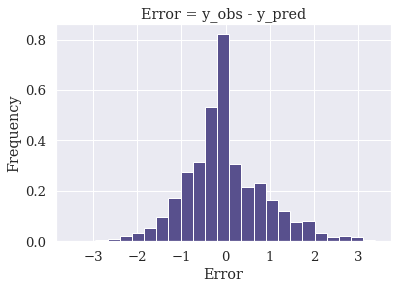}
  \caption{}
  \label{fig:a_sub2}
\end{subfigure}

\centering
\begin{subfigure}{.4\textwidth}
  \centering
  \includegraphics[width=1\linewidth]{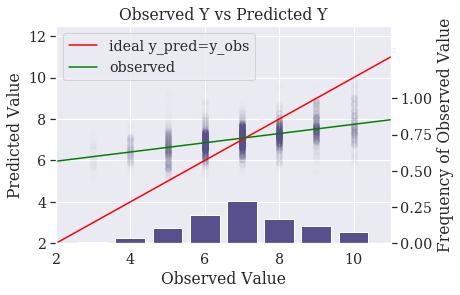}
  \caption{}
  \label{fig:scatter_b_sub1}
\end{subfigure}%
\begin{subfigure}{.4\textwidth}
  \centering
  \includegraphics[width=1\linewidth]{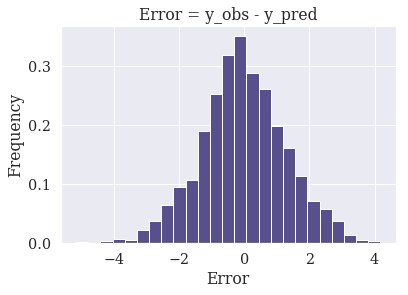}
  \caption{}
  \label{fig:b_sub2}
\end{subfigure}

\centering
\begin{subfigure}{.4\textwidth}
  \centering
  \includegraphics[width=1\linewidth]{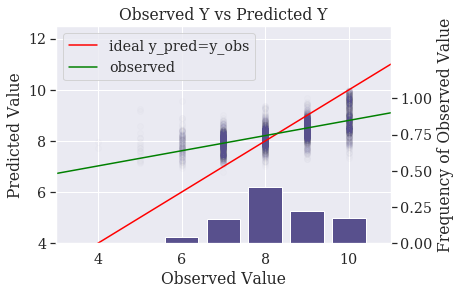}
  \caption{}
  \label{fig:scatter_c_sub1}
\end{subfigure}%
\begin{subfigure}{.4\textwidth}
  \centering
  \includegraphics[width=1\linewidth]{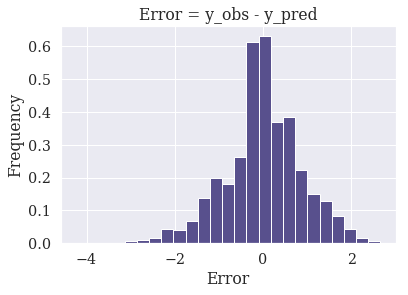}
  \caption{}
  \label{fig:c_sub2}
\end{subfigure}
\caption{\textbf{References} $A$ (a, b), $B$ (c,d) and $C$ (e,f): Learning, Regression}
\label{fig:ml_capc1}
\end{figure}

Figures \ref{fig:scatter_a_sub1}, \ref{fig:scatter_b_sub1} and \ref{fig:scatter_c_sub1}, plot the predictions of $Y$ (number of facilities kept open) versus the true values on the test data.
The ideal case of perfect predictions $\hat{Y} = Y$, is denoted by the red line.
The green line is a simple linear regression of the form $\widehat Y = \hat m Y + \hat b$ that graphically renders the actual relation between $Y$ and $\widehat Y$.

We remark that the data are highly concentrated around some central values. 
The tails of the distribution are scarcely represented in our training data, and the learning function cannot generalize well, as one can see more prominently with instance $A$, where most of the data take values 2, 3 or 4.
\subsection{Optimization Results}\label{sec:results}
We report the performance of our method on the three references $A, B$ and $C$, generating three batches $D_A, D_B$ and $D_C$ of $50$ new instances, unseen during training.
We apply the same procedure followed to generate the learning data (Section \ref{sec:gendata}).

To measure the effect of the ML bound with respect to resolution time, each instance in $D_A$ (but also $D_B$ and $D_C$) is optimized twice (with and without bound) with time limits $t_{lim} \in \{5, 10, 30,60, 120\}$ seconds.
At the end of the time-limited run, we record the objective value and the effective solving time (some instances are optimized before the limit is reached).

Our ML constraint cannot make the solution infeasible, as the total number of active facilities is not bounded, but it could cut off the optimal solution.
For instance, one could predict all of the reference facilities to be open ($\hat \gamma = 1$) while none of them should be.
We run the optimization with 3,600 seconds of time limit, to measure the difference in objective values after a long run, with and without the additional constraint.

\subsubsection{Experimental setup}
For each iteration of our experiments, we proceed as follows.
A new instance is generated as detailed in Section \ref{sec:gendata}.
The instance is then solved once for each time limit without the additional constraint.
The objectives and effective times are recorded. 
The ML algorithm is run (see Section \ref{sec:learnbound}) and constraint \eqref{eq:constr} is introduced in the model's formulation.
The constrained model is then solved once for each time limit.
To verify that the optimal solution was preserved by the introduced bound constraint, we run both models with and without constraint \eqref{eq:constr} for a full hour to check if any remarkable deviation has been caused by the ML.

When using a highly-optimized MIP solver like CPLEX, it is sometimes hard to evaluate the impact of model modifications (imposed by the user) because the complex algorithmic components within the solver might mitigate the effect of those modifications. For this reason, it is common practice to deactivate, for testing purposes, some of those CPLEX algorithmic ingredients. 

In our special case, we are interested in evaluating the effect of the ML-predicted constraint \eqref{eq:constr} as a way to fast dive toward good solutions of a modified FLP instance with respect to the solution of a reference one. To achieve this goal, we have found useful to run CPLEX without the \textit{presolve} feature and we first report the results of this configuration. We nonetheless performed all experiments and report the results with CPLEX \textit{default}, i.e.,  with \textit{presolve} activated. We show that the overall message provided by the experiments is in both cases very similar.

Finally, given the interest, especially in the academic community, of using non-commercial MIP solvers, we report in Appendix \ref{annex:scip} the results obtained by replacing CPLEX with the most widely adopted of those solvers, namely SCIP (\textit{default} version) \citep{GleixnerEtal2017OO}.

\subsubsection{Experimental results}

Within $3,600$ seconds, most of the runs attain an optimality gap of the order of $0.05\%$ or manage to reach optimality.
The bias introduced by the learned bound is on average around $0.03\%$ and at most of the order of $0.5\%$.
That is, at the $3,600$ seconds time limit, the objective value of the ML-bounded run is on average $0.03\%$ higher than the objective for the regular run, but never greater than $0.5\%$.

\begin{figure}
\centering
\begin{subfigure}{.4\textwidth}
  \centering
  \includegraphics[width=1\linewidth]{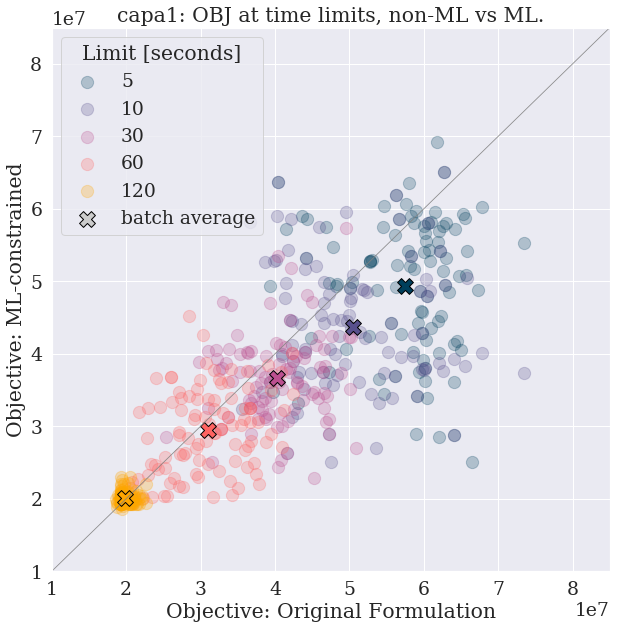}
  \caption{}
  \label{fig:a_objobj}
\end{subfigure}%
\begin{subfigure}{.4\textwidth}
  \centering
  \includegraphics[width=1\linewidth]{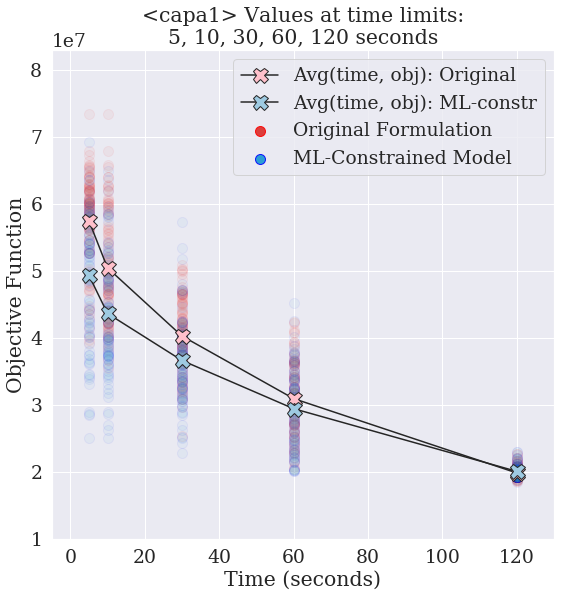}
  \caption{}
  \label{fig:a_timeobj}
\end{subfigure}
 
\centering
\begin{subfigure}{.4\textwidth}
  \centering
  \includegraphics[width=1\linewidth]{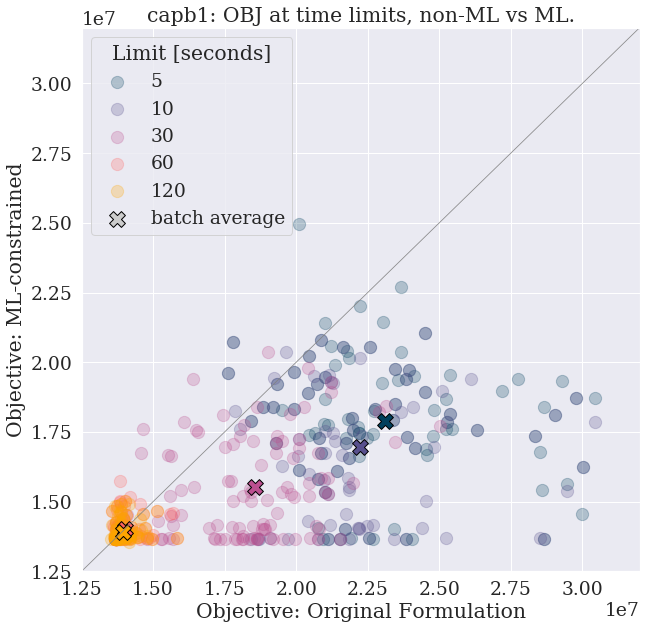}
  \caption{ }
  \label{fig:b_objobj}
\end{subfigure}%
\begin{subfigure}{.4\textwidth}
  \centering
  \includegraphics[width=1\linewidth]{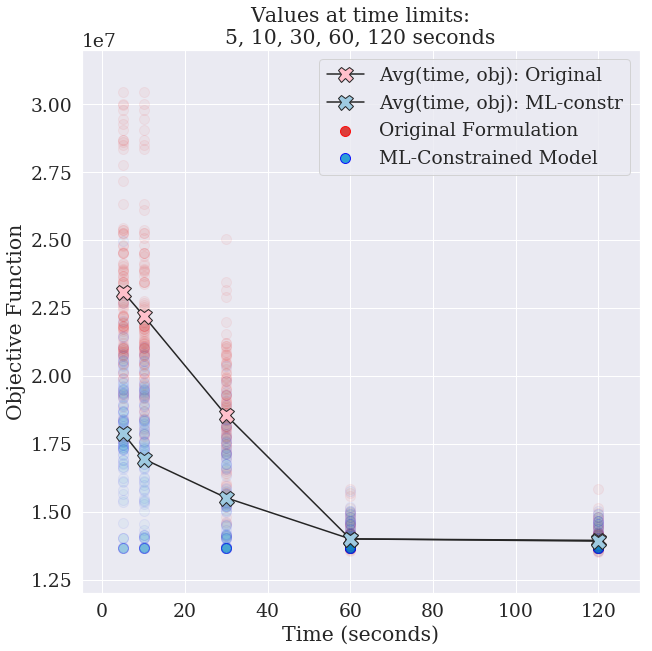}
  \caption{} 
  \label{fig:b_timeobj}
\end{subfigure}

\centering
\begin{subfigure}{.4\textwidth}
  \centering
  \includegraphics[width=1\linewidth]{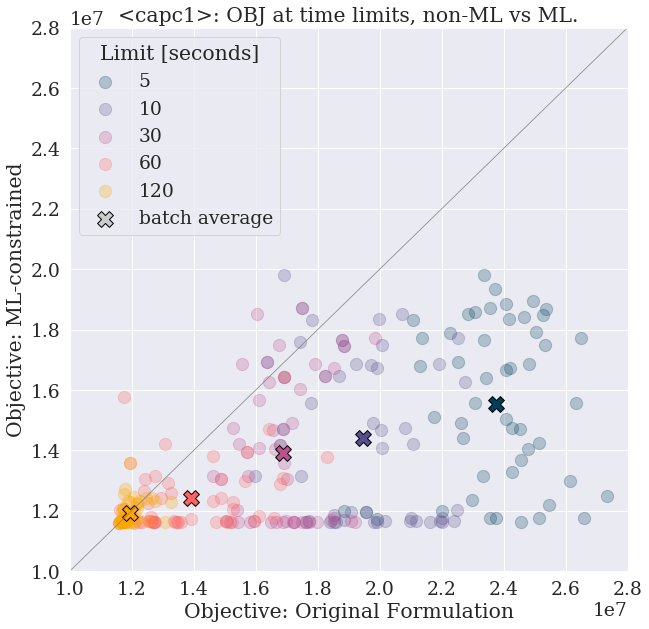}
  \caption{} 
  \label{fig:c_objobj}
\end{subfigure}%
\begin{subfigure}{.4\textwidth}
  \centering
  \includegraphics[width=1\linewidth]{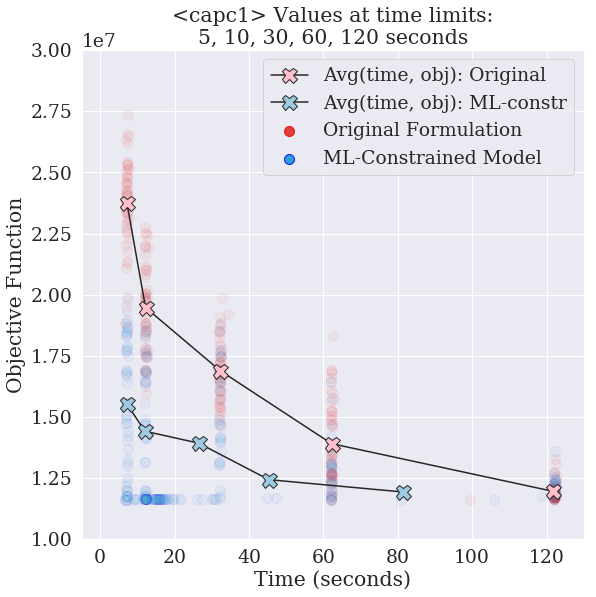}
  \caption{} 
  \label{fig:c_timeobj}
\end{subfigure}
\caption{\textbf{References} $A$ (a, b), $B$ (c,d) and $C$ (e,f): optimization results}
\end{figure}

Figures \ref{fig:a_objobj}, \ref{fig:b_objobj} and \ref{fig:c_objobj} compare the values of the objective functions of the constrained and unconstrained resolutions at the time limit.
Each point represents the same new, unseen, instance optimized twice: once with the ML bound and once without.
The color coding helps determining the time limit imposed on this pair of runs.
On the $x$ axis one can see the objective value without ML, on the $y$ axis the value with the ML bound.
Points on the lower left corner correspond to runs with longer time limits, where the difference between the two approaches becomes more negligible.
The points below the line are runs for which our method outperformed the original formulation.

In Figures \ref{fig:a_timeobj}, \ref{fig:b_timeobj} and \ref{fig:c_timeobj}, the same information is presented focusing on the temporal dimension.
On the $x$ axis is the effective resolution time, which can be inferior to the time limit.
Along the $y$ axis is the objective value at the time limit.
The points aligned at the bottom of Figures \ref{fig:b_timeobj} and \ref{fig:c_timeobj} are those for which the optimal solution was found before the time limit expired.
For short resolution times, the ML-constrained formulation yields better objective function values than the original formulation without cutting the optimal solution which is eventually found given more computational time.
On average, the blue points representing the objective of the ML-bounded model lie below the red ones of the unbounded model.
In case $C$, our method allows for superior performances than those obtained simply via the solver.
In case $A$ we are still on average improving the performances or at least keeping in line with the unbounded model.
Here, the amplitude of the perturbations imply completely different solutions, making transfer difficult.
We remark however that, on average, our approach avoids cutting the optimal solution and thus prevents \textit{negative} transfer.

\begin{figure}
\centering
\begin{subfigure}{.4\textwidth}
  \centering
  \includegraphics[width=1\linewidth]{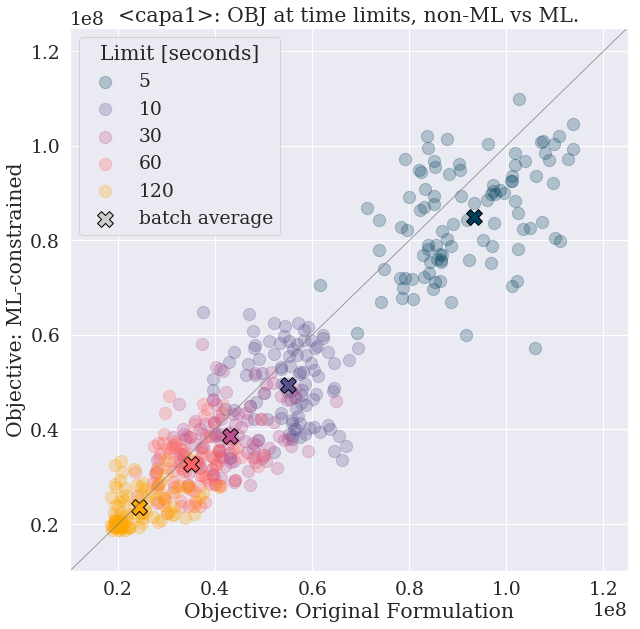}
  \caption{}
\end{subfigure}%
\begin{subfigure}{.4\textwidth}
  \centering
  \includegraphics[width=1\linewidth]{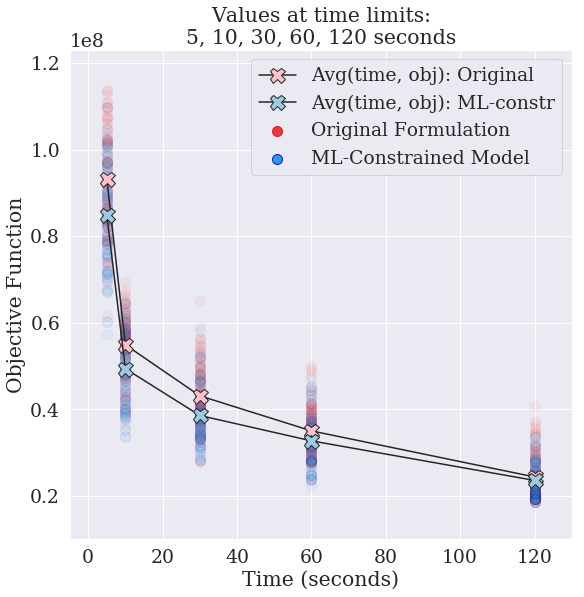}
  \caption{}
\end{subfigure}

\centering
\begin{subfigure}{.4\textwidth}
  \centering
  \includegraphics[width=1\linewidth]{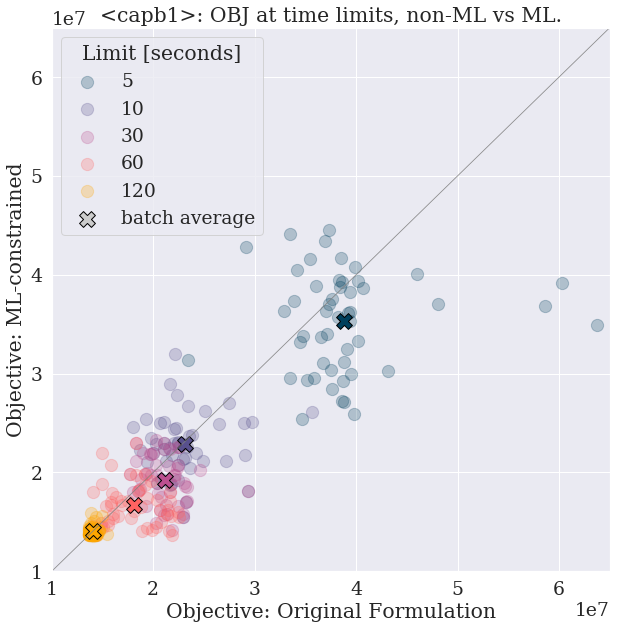}
  \caption{}
\end{subfigure}%
\begin{subfigure}{.4\textwidth}
  \centering
  \includegraphics[width=1\linewidth]{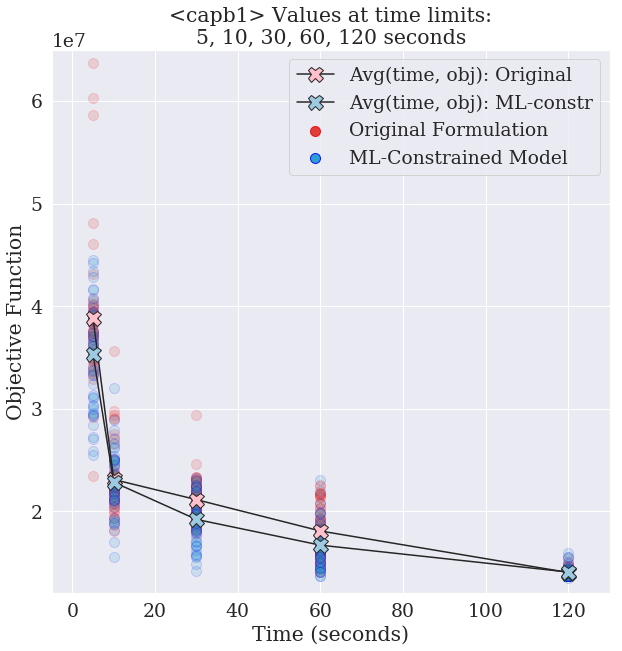}
  \caption{}
\end{subfigure}

\centering
\begin{subfigure}{.4\textwidth}
  \centering
  \includegraphics[width=1\linewidth]{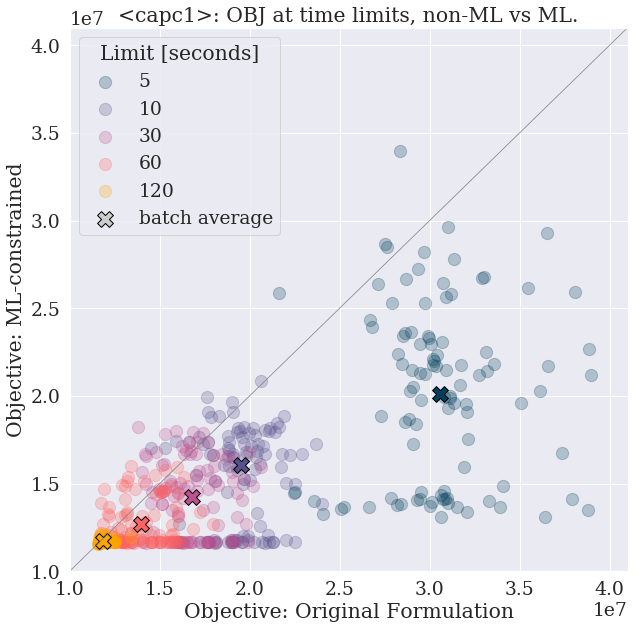}
  \caption{}
\end{subfigure}%
\begin{subfigure}{.4\textwidth}
  \centering
  \includegraphics[width=1\linewidth]{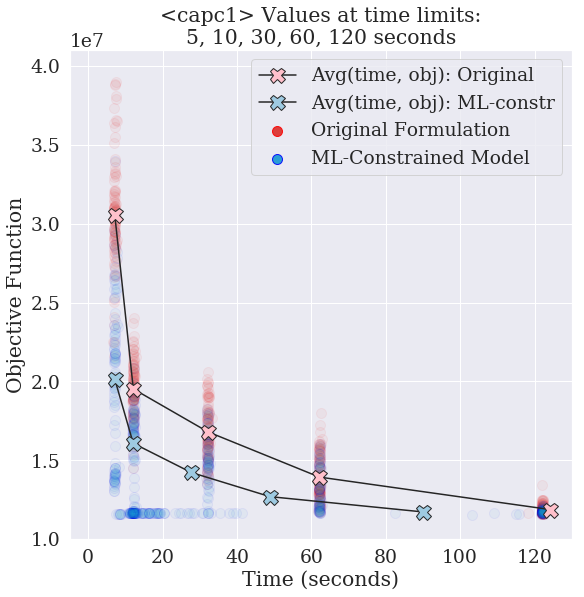}
  \caption{}
\end{subfigure}
\caption{\textbf{References} $A$ (a, b), $B$ (c,d) and $C$ (e,f): optimization results with \textit{presolve}}
\label{fig:preso}
\end{figure}

As noted above, we operated with the \textit{presolve} feature switched off.
For completeness, we run all of the experiments a second time with the \textit{presolve} feature activated (Figure \ref{fig:preso}).
The results reflect what was observed previously, although the presolve heuristics of CPLEX manage to reduce in part the effect of our approach in case $B$ (Figure \ref{fig:preso} c and d).

\section{Conclusions and Perspectives}
\label{sec:conclusion}

We have shown that the existing or simulated data of a recurring optimization problem can be exploited to gain insight into the optimization process.
We have used these data to fit a binary classifier and a regressor.
We predict whether a new instance, derived from a reference instance, will share all or a part of its optimal solution with the solution of the reference instance.
This piece of information, translated into an additional constraint, is given to a solver and allows to dive faster, on average, towards a good solution.
Moreover, we empirically illustrated that this additional constraint preserves the optimal solution and thus prevents negative transfer.

Building upon these results, we plan on extending the reach of our approach by experimenting with other families of problems, where, potentially, different types of constraints should be considered. This will be the way of fully demonstrating the wide applicability of the approach to both solving repetition (with different data) and re-optimization.

\section*{Acknowledgements}
The authors wish to acknowledge the support of the Canada Excellence Research Chair in Data Science for Real-Time Decision-Making at Polytechnique Montr\'{e}al and the ISAE-SUPAERO foundation.

\bibliography{facility} 

\begin{thebibliography}{}

\bibitem[Alvarez et~al., 2017]{alvarez2017machine}
Alvarez, A.~M., Louveaux, Q., and Wehenkel, L. (2017).
\newblock A machine learning-based approximation of strong branching.
\newblock {\em INFORMS Journal on Computing}, 29(1):185--195.

\bibitem[Basso et~al., 2018]{Basso2018}
Basso, S., Ceselli, A., and Tettamanzi, A. (2018).
\newblock Random sampling and machine learning to understand good
  decompositions.
\newblock {\em Annals of Operations Research}.

\bibitem[Bello et~al., 2017]{bello2017neural}
Bello, I., Zoph, B., Vasudevan, V., and Le, Q.~V. (2017).
\newblock Neural optimizer search with reinforcement learning.
\newblock In {\em International Conference on Machine Learning}, pages
  459--468.

\bibitem[Bengio et~al., 2018]{BLP18}
Bengio, Y., Lodi, A., and Prouvost, A. (2018).
\newblock Machine learning for combinatorial optimization: a methodological
  tour d'horizon.
\newblock {\em Preprint arXiv:1811.06128}.

\bibitem[Bonami et~al., 2018]{bonami2018learning}
Bonami, P., Lodi, A., and Zarpellon, G. (2018).
\newblock Learning a classification of mixed-integer quadratic programming
  problems.
\newblock In van Hoeve, W.-J., editor, {\em Integration of Constraint
  Programming, Artificial Intelligence, and Operations Research}, pages
  595--604, Cham. Springer International Publishing.

\bibitem[Cornu{\'e}jols et~al., 1983]{cornuejols1983uncapacitated}
Cornu{\'e}jols, G., Nemhauser, G.~L., and Wolsey, L.~A. (1983).
\newblock The uncapacitated facility location problem.
\newblock Technical report, Carnegie-mellon univ pittsburgh pa management
  sciences research group.

\bibitem[Deudon et~al., 2018]{deudon2018learning}
Deudon, M., Cournut, P., Lacoste, A., Adulyasak, Y., and Rousseau, L.-M.
  (2018).
\newblock Learning heuristics for the tsp by policy gradient.
\newblock In {\em International Conference on the Integration of Constraint
  Programming, Artificial Intelligence, and Operations Research}, pages
  170--181. Springer.

\bibitem[Di~Liberto et~al., 2016]{Liberto2016DASH}
Di~Liberto, G., Kadioglu, S., Leo, K., and Malitsky, Y. (2016).
\newblock Dash: Dynamic approach for switching heuristics.
\newblock {\em European Journal of Operational Research}, 248(3):943 -- 953.

\bibitem[Fischetti and Fraccaro, 2018]{fischetti2018machine}
Fischetti, M. and Fraccaro, M. (2018).
\newblock Machine learning meets mathematical optimization to predict the
  optimal production of offshore wind parks.
\newblock {\em Computers and Operations Research}.

\bibitem[Fischetti and Lodi, 2003]{localbranching}
Fischetti, M. and Lodi, A. (2003).
\newblock Local branching.
\newblock {\em Mathematical Programming}, 98:23--47.

\bibitem[Garey and Johnson, 1990]{Garey1979ComputersAI}
Garey, M.~R. and Johnson, D.~S. (1990).
\newblock {\em Computers and Intractability; A Guide to the Theory of
  NP-Completeness}.
\newblock W. H. Freeman \& Co., New York, NY, USA.

\bibitem[Geurts et~al., 2006]{geurts2006extremely}
Geurts, P., Ernst, D., and Wehenkel, L. (2006).
\newblock Extremely randomized trees.
\newblock {\em Machine learning}, 63(1):3--42.

\bibitem[Gleixner et~al., 2017]{GleixnerEtal2017OO}
Gleixner, A., Eifler, L., Gally, T., Gamrath, G., Gemander, P., Gottwald,
  R.~L., Hendel, G., Hojny, C., Koch, T., Miltenberger, M., M{\"u}ller, B.,
  Pfetsch, M.~E., Puchert, C., Rehfeldt, D., Schl{\"o}sser, F., Serrano, F.,
  Shinano, Y., Viernickel, J.~M., Vigerske, S., Weninger, D., Witt, J.~T., and
  Witzig, J. (2017).
\newblock {The SCIP Optimization Suite 5.0}.
\newblock Technical report, Optimization Online.

\bibitem[Goodfellow et~al., 2016]{Goodfellow-et-al-2016}
Goodfellow, I., Bengio, Y., and Courville, A. (2016).
\newblock {\em Deep Learning}.
\newblock MIT Press.
\newblock \url{http://www.deeplearningbook.org}.

\bibitem[Guastaroba and Speranza, 2014]{GUASTAROBA2014_heuristic}
Guastaroba, G. and Speranza, M. (2014).
\newblock A heuristic for bilp problems: The single source capacitated facility
  location problem.
\newblock {\em European Journal of Operational Research}, 238(2):438 -- 450.

\bibitem[Guastaroba and Speranza, 2012]{guastaroba2012kernel}
Guastaroba, G. and Speranza, M.~G. (2012).
\newblock Kernel search for the capacitated facility location problem.
\newblock {\em Journal of Heuristics}, 18(6):877--917.

\bibitem[Hastie et~al., 2009]{hastie2009unsupervised}
Hastie, T., Tibshirani, R., and Friedman, J. (2009).
\newblock The elements of statistical learning: data mining, inference, and
  prediction, springer series in statistics.

\bibitem[He et~al., 2014]{he2014learning}
He, H., Daume~III, H., and Eisner, J.~M. (2014).
\newblock Learning to search in branch and bound algorithms.
\newblock In {\em Advances in Neural Information Processing Systems}, pages
  3293--3301.

\bibitem[IBM, 2018]{cplex_solver}
IBM (2018).
\newblock {\em IBM ILOG CPLEX Optimizers 12.8.0}.

\bibitem[Khalil et~al., 2017]{khalil2017learning}
Khalil, E.~B., Dai, H., Zhang, Y., Dilkina, B., and Song, L. (2017).
\newblock Learning combinatorial optimization algorithms over graphs.
\newblock In {\em Advances in Neural Information Processing Systems}, pages
  6348--6358.

\bibitem[Khalil et~al., 2016]{KhaBodSonNemDil16}
Khalil, E.~B., Le~Bodic, P., Song, L., Nemhauser, G., and Dilkina, B. (2016).
\newblock Learning to branch in mixed integer programming.
\newblock In {\em Proceedings of the 30th AAAI Conference on Artificial
  Intelligence}.

\bibitem[Klose and Drexl, 2005]{klose2005facility}
Klose, A. and Drexl, A. (2005).
\newblock Facility location models for distribution system design.
\newblock {\em European Journal of Operational Research}, 162(1):4 -- 29.
\newblock Logistics: From Theory to Application.

\bibitem[Krizhevsky et~al., 2012]{krizhevsky2012imagenet}
Krizhevsky, A., Sutskever, I., and Hinton, G.~E. (2012).
\newblock Imagenet classification with deep convolutional neural networks.
\newblock In {\em Advances in neural information processing systems}, pages
  1097--1105.

\bibitem[Kruber et~al., 2017]{kruber2017learning}
Kruber, M., L{\"u}bbecke, M.~E., and Parmentier, A. (2017).
\newblock Learning when to use a decomposition.
\newblock In Salvagnin, D. and Lombardi, M., editors, {\em Integration of AI
  and OR Techniques in Constraint Programming}, pages 202--210, Cham. Springer
  International Publishing.

\bibitem[Larsen et~al., 2018]{larsen2018predicting}
Larsen, E., Lachapelle, S., Bengio, Y., Frejinger, E., Lacoste-Julien, S., and
  Lodi, A. (2018).
\newblock Predicting solution summaries to integer linear programs under
  imperfect information with machine learning.
\newblock {\em arXiv preprint arXiv:1807.11876}.

\bibitem[Lodi and Zarpellon, 2017]{lodi2017learning}
Lodi, A. and Zarpellon, G. (2017).
\newblock On learning and branching: a survey.
\newblock {\em TOP}, pages 1--30.

\bibitem[Mnih et~al., 2015]{mnih2015human}
Mnih, V., Kavukcuoglu, K., Silver, D., Rusu, A.~A., Veness, J., Bellemare,
  M.~G., Graves, A., Riedmiller, M., Fidjeland, A.~K., Ostrovski, G., et~al.
  (2015).
\newblock Human-level control through deep reinforcement learning.
\newblock {\em Nature}, 518(7540):529.

\bibitem[Xavier et~al., 2019]{Ahmed2019}
Xavier, A., Qiu, F., and Ahmed, S. (2019).
\newblock Learning to solve large-scale security-constrained unit commitment
  problems.
\newblock Technical Report 1902.01697, arXiv.

\end{thebibliography}
\bibliographystyle{apalike}

\appendix
\section{Results with SCIP} 
\label{annex:scip}

\begin{figure}
\centering
\begin{subfigure}{.4\textwidth}
  \centering
  \includegraphics[width=1\linewidth]{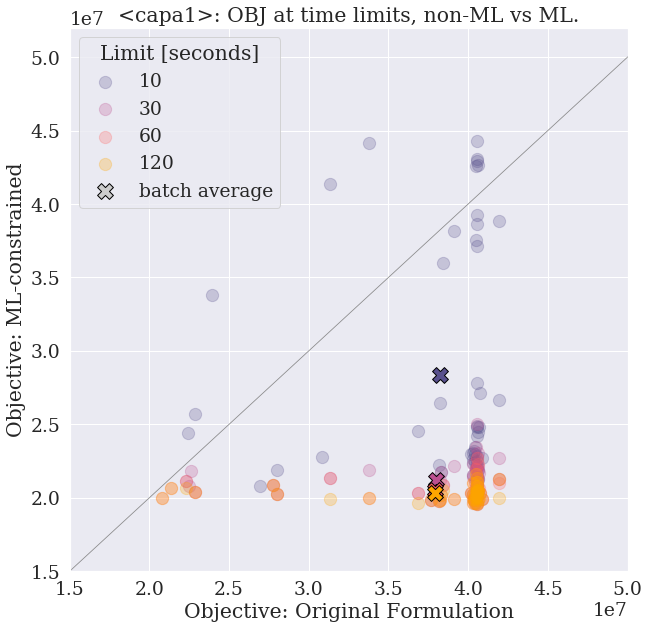}
  \caption{}
\end{subfigure}%
\begin{subfigure}{.4\textwidth}
  \centering
  \includegraphics[width=1\linewidth]{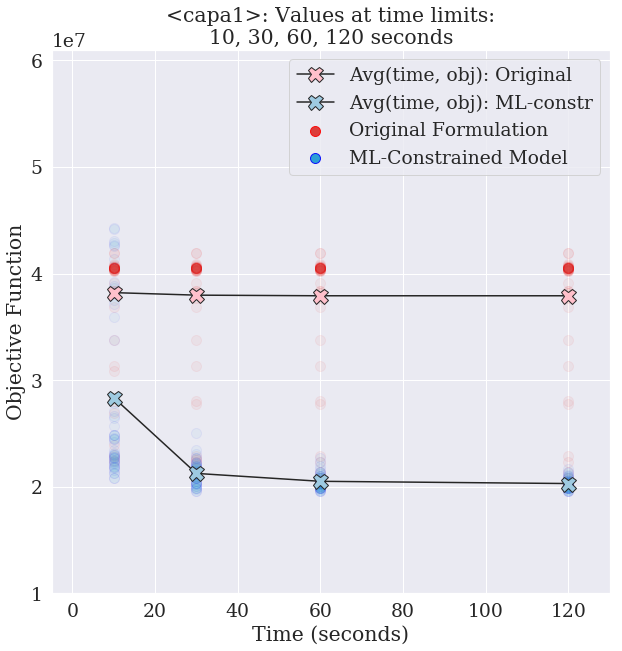}
  \caption{}
\end{subfigure}

\centering
\begin{subfigure}{.4\textwidth}
  \centering
  \includegraphics[width=1\linewidth]{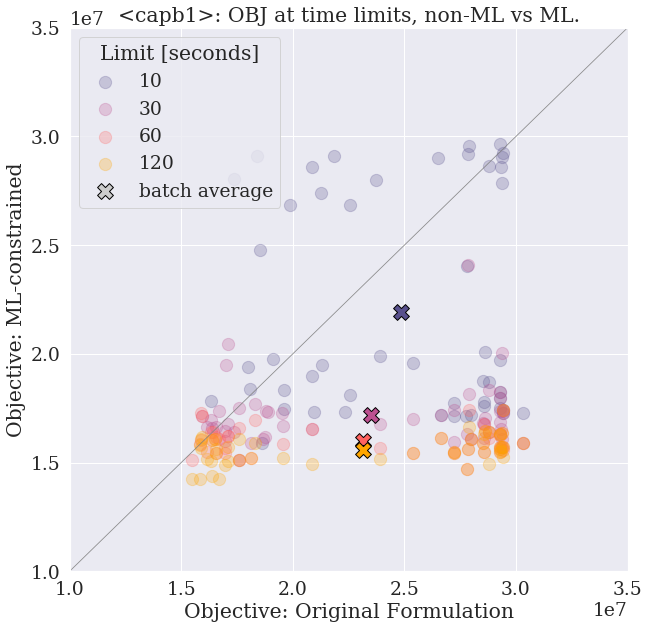}
  \caption{}
\end{subfigure}%
\begin{subfigure}{.4\textwidth}
  \centering
  \includegraphics[width=1\linewidth]{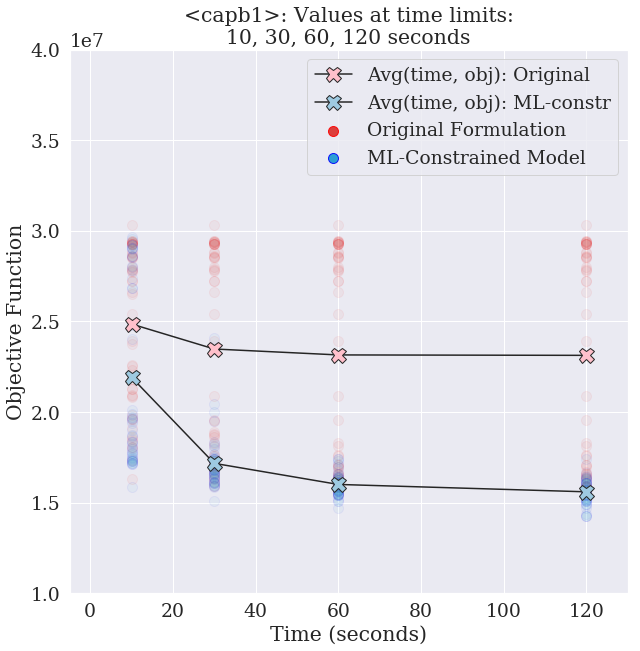}
  \caption{}
\end{subfigure}
\centering
\begin{subfigure}{.4\textwidth}
  \centering
  \includegraphics[width=1\linewidth]{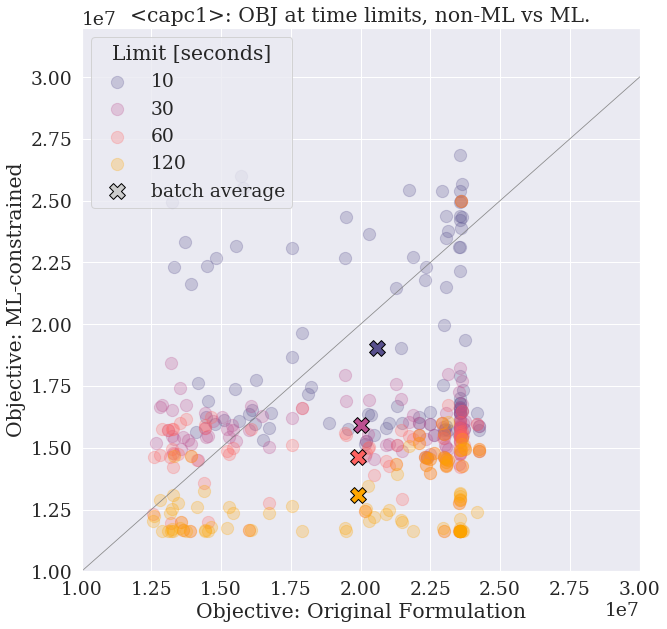}
  \caption{}
\end{subfigure}%
\begin{subfigure}{.4\textwidth}
  \centering
  \includegraphics[width=1\linewidth]{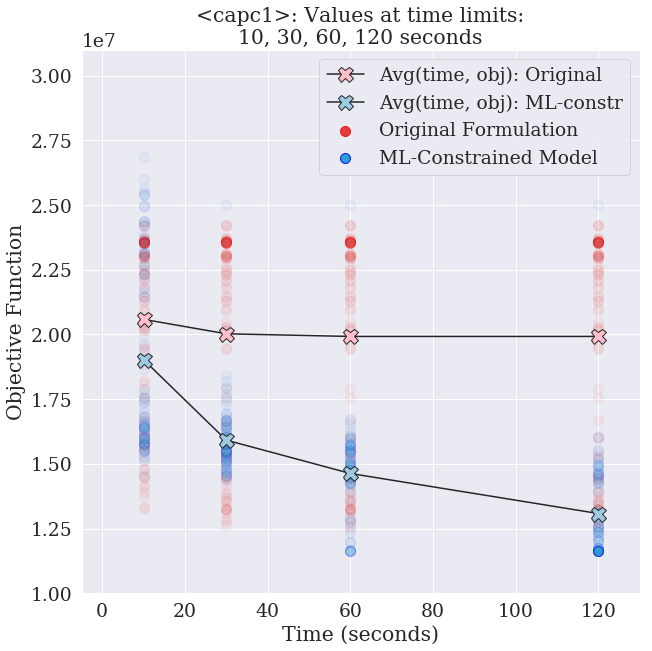}
  \caption{}
\end{subfigure}
\caption{\textbf{References} $A$ (a, b), $B$ (c,d) and $C$ (e,f): optimization results with the SCIP solver \label{fig:SCIP}}
\end{figure}

The MIP solver SCIP, whose source code is accessible for academic purposes, has a considerable adoption in academic research, which is why it is interesting to test our approach in conjunction with such a solver.
The experiments reported in Figure \ref{fig:SCIP} are the same as those in Section \ref{sec:results}, with the exception of the time limits of $5$ seconds; for such a short time limit, most of the runs did not yield a feasible solution and we did not report it.
What observed previously with CPLEX is true here as well, although the effect of constraint \eqref{eq:constr} are amplified.
The additional bound constraint allows the solver to dive much faster towards a good solution.
Even at the 120 seconds mark, the gap between the two objective values is still considerable.

\end{document}